\RequirePackage{fix-cm}
\documentclass[smallextended]{svjour3}       % onecolumn (second format)
\smartqed  % flush right qed marks, e.g. at end of proof
\usepackage{graphicx}
\graphicspath{{figs/}}
\usepackage{subfigure}
\usepackage{amsmath}
\usepackage{color}
\usepackage{amsfonts,amssymb}
\usepackage{hyperref}
\usepackage{array}  
\usepackage{booktabs} %调整表格线与上下内容的间隔  
\usepackage{multirow}
\usepackage{algpseudocode}
\usepackage{amsmath}
\usepackage[ruled]{algorithm2e}

\hypersetup{hypertex=true,
 colorlinks=true,
linkcolor=blue,
anchorcolor=blue,
citecolor=blue}
\allowdisplaybreaks[4]

\begin{document}

\title{A gradient method exploiting the two dimensional quadratic termination property%\thanks{Grants or other notes
}

\titlerunning{A gradient method exploiting the 2D quadratic termination property}        % if too long for running head

\author{Xinrui Li         \and
        Yakui Huang %etc.
}


\institute{Xinrui Li \at
             School of Sciences, Hebei University of Technology, Tianjin 300401, China 
%             \\
%              \email{xinruili1020@163.com}           %  \\
           \and
           Yakui Huang \at
              Institute of Mathematics, Hebei University of Technology, Tianjin 300401, China\\
             \email{ huangyakui2006@gmail.com}
}

\date{Received: date / Accepted: date}

\maketitle

\begin{abstract}
The quadratic termination property is important to the efficiency of gradient methods. We consider equipping a family of gradient methods, where the stepsize is given by the ratio of two norms, with two dimensional quadratic termination. Such a desired property is achieved by cooperating with a new stepsize which is derived by maximizing the stepsize of the considered family in the next iteration. It is proved that each method in the family will asymptotically alternate in a two dimensional subspace spanned by the eigenvectors corresponding to the largest and smallest eigenvalues. Based on this asymptotic behavior, we show that the new stepsize converges to the reciprocal of the largest eigenvalue of the Hessian. Furthermore, by adaptively taking the long Barzilai--Borwein  stepsize and  reusing the new stepsize with retard, we propose an efficient gradient method for unconstrained quadratic optimization. We prove that the new method is $R$-linearly convergent with a rate of $1-1/\kappa$, where $\kappa$ is the condition number of Hessian. Numerical experiments show the efficiency of our proposed method.
%A new property is presented to establish a uniform $R$-linear convergence result for gradient methods, which implies the proposed method and many others are $R$-linearly convergent. Under proper conditions, we prove the $R$-factor can be bounded above by $1-1/\kappa$, where $\kappa$ is the condition number of Hessian. 

\end{abstract}

\section{Introduction}
\label{intro}
The gradient method is well-known for minimizing a smooth  function
$f(x): \mathbb{R}^{n} \rightarrow \mathbb{R}$, which updates the iterates by
\begin{equation}\label{aa}
	x_{k+1}=x_{k}-\alpha_{k}g_{k},
\end{equation}
where $g_{k} = \nabla f\left(x_{k}\right)$ and the stepsize $\alpha_{k}>0$ depends on the method under consideration. The classic steepest descent (SD) \cite{cauchy1847methode} and minimal gradient (MG)  \cite{dai2003altermin} methods  determine $\alpha_{k}$ by minimizing $f\left(x_{k}-\alpha g_{k}\right)$ and $\left\|g\left(x_{k}-\alpha g_{k}\right)\right\|_2$, respectively. Theoretically, the SD and MG methods will asymptotically perform zigzags between two directions, which often yield poor performance in many probelms  \cite{akaike1959successive,forsythe1968asymptotic,huang2019asymptotic,nocedal2002behavior,zou2021fast}.

 In 1988, Barzilai and Borwein \cite{barzilai1988two} proposed the following two novel stepsizes  from the view of quasi-Newton methods,
 \begin{equation}\label{bb1}
 \alpha_{k}^{B B 1}=\arg \min _{\alpha \in \mathbb{R}}\left\|\alpha^{-1} s_{k-1}-y_{k-1}\right\|_{2}=\frac{s_{k-1}^{T} s_{k-1}}{s_{k-1}^{T} y_{k-1}}
 \end{equation}
and
 \begin{equation}\label{bb2}
\alpha_{k}^{B B 2}=\arg \min _{\alpha \in \mathbb{R}}\left\|s_{k-1}-\alpha y_{k-1}\right\|_{2}=\frac{s_{k-1}^{T} y_{k-1}}{y_{k-1}^{T} y_{k-1}},
\end{equation}
where $s_{k-1}=x_{k}-x_{k-1}$ and $y_{k-1}=g_{k}-g_{k-1}$. Apparently, $\alpha_{k}^{BB1} \geq \alpha_{k}^{BB2}$ follows from the Cauchy--Schwartz
inequality if $s_{k-1}^{T}y_{k-1} > 0$. When  $f(x)$ is a quadratic function
\begin{equation}\label{1a}
	f\left(x\right)=\frac{1}{2}x^{T}Ax-b^{T}x,
\end{equation}
where $A\in \mathbb{R}^{n\times n}$ is a  real  symmetric positive definite  matrix
and $b\in \mathbb{R}^{n}$, $\alpha_{k}^{B B 1}$ and $\alpha_{k}^{B B 2}$ can be regarded as the SD and MG stepsizes with retard, respectively, i.e.
\begin{equation*}
	\alpha_{k}^{B B 1}=\frac{g_{k-1}^{T} g_{k-1}}{g_{k-1}^{T} A g_{k-1}}=\alpha_{k-1}^{SD} \quad \text { and } \quad \alpha_{k}^{B B 2}=\frac{g_{k-1}^{T} A g_{k-1}}{g_{k-1}^{T} A^{2} g_{k-1}}=\alpha_{k-1}^{MG}.
\end{equation*}
Barzilai and Borwein \cite{barzilai1988two} proved that the BB method is $R$-superlinearly convergent for the two dimensional strictly convex quadratic function.
It has been shown that the BB method is globally convergent \cite{raydan1993barzilai} with $R$-linear rate \cite{dai2002r} for any  dimensional cases.  Although the BB method is nonmonotone,  extensive numerical experimental results indicate that it performs much better than the SD method \cite{fletcher2005barzilai,raydan1997barzilai,yuan2008step}. See \cite{birgin2000nonmonotone,birgin2014spectral,dai2019family,di2018steplength,grippo1986nonmonotone,huang2016smoothing,huang2015quadratic,jiang2013feasible}  for more BB-like methods.

In \cite{yuan2006new}, Yuan derived a new stepsize,
%\begin{equation}\label{yuanstepsize}
%	\alpha_{k}^{\mathrm{Y}}=2\left(\sqrt{\left(\frac{1}{\alpha_{k-1}^{\mathrm{SD}}}-\frac{1}{\alpha_{k}^{\mathrm{SD}}}\right)^{2}+\frac{4\left\|g_{k}\right\|_{2}^{2}}{\left\|s_{k-1}\right\|_{2}^{2}}}+\frac{1}{\alpha_{k-1}^{\mathrm{SD}}}+\frac{1}{\alpha_{k}^{\mathrm{SD}}}\right)^{-1},
%\end{equation}
which together with the SD method  produces the
minimizer of a two dimensional strictly convex quadratic function in three iterations. In what follows, if a
method can give the exact minimizer of a two dimensional convex quadratic
function within finite iterations, we call it has the property of two dimensional
quadratic termination.  Based on the following variant of the Yuan stepsize,
 \begin{equation}\label{dystepsize}
% 	\alpha_{k}^{\mathrm{DY}}=2\left(\sqrt{\left(\frac{1}{\alpha_{k-1}^{\mathrm{SD}}}-\frac{1}{\alpha_{k}^{\mathrm{SD}}}\right)^{2}+\frac{4\left\|g_{k}\right\|_2^{2}}{\left(\alpha_{k-1}^{\mathrm{SD}}\right)^{2}\left\|g_{k-1}\right\|_2^{2}}}+\frac{1}{\alpha_{k-1}^{\mathrm{SD}}}+\frac{1}{\alpha_{k}^{\mathrm{SD}}}\right)^{-1},
 	\alpha_{k}^{\mathrm{DY}}=\frac{2}{\sqrt{\left(\frac{1}{\alpha_{k-1}^{\mathrm{SD}}}-\frac{1}{\alpha_{k}^{\mathrm{SD}}}\right)^{2}+\frac{4\left\|g_{k}\right\|_2^{2}}{\left(\alpha_{k-1}^{\mathrm{SD}}\right)^{2}\left\|g_{k-1}\right\|_2^{2}}}+\frac{1}{\alpha_{k-1}^{\mathrm{SD}}}+\frac{1}{\alpha_{k}^{\mathrm{SD}}}},
 \end{equation}
Dai and Yuan \cite{dai2005analysis} suggested the so-called Dai--Yuan (DY) gradient method with
\begin{equation}
	\alpha_{k} =\left\{\begin{array}{ll}
		\alpha_{k}^{S D}, & \text { if } \bmod (k, 4)<2; \\
		\alpha_{k}^{D Y}, & \text { otherwise. }
	\end{array}\right.
\end{equation}
Clearly, $\alpha_{k}^{DY}\leq \min\{\alpha_{k}^{SD},\alpha_{k-1}^{SD}\}$, which implies that the DY method is monotone. Interestingly, the DY method can even outperform the nonmonotone BB method. Recently, Huang et al. \cite{huang2019asymptotic} derived a new stepsize, say $\alpha_{k}^{H}$,  such that the gradient method
 $$\alpha_{k} = \frac{g_{k}^{T} \psi(A) g_{k}}{g_{k}^{T} \psi(A) A g_{k}}$$ 
together with $\alpha_{k}^{H}$ achieves the two dimensional quadratic termination, where $\psi$ is a real analytic function  on $[\lambda_{1},\lambda_{n}]$ and  can be expressed by Laurent series 
\begin{equation}\label{psi}
	\psi(z)=\sum_{k=-\infty}^{\infty} d_{k} z^{k},~d_{k} \in \mathbb{R},
\end{equation}
such that $0<\sum_{k=-\infty}^{\infty} d_{k} z^{k}<+\infty$ for all $z \in\left[\lambda_{1}, \lambda_{n}\right]$. Here, $\lambda_{1}$ and $\lambda_{n}$ are the smallest and largest eigenvalues of $A$, respectively. Furthermore,  $\alpha_{k}^{H}$ reduces to $\alpha_{k}^{DY}$ when $\psi(A)=I$. The property of two dimensional quadratic termination has shown great potential in improving performances of gradient methods, see \cite{huang2020equipping,huang2019asymptotic,sun2020new} for example.

%asymptotic optimal gradient (Dai--Yang) method  designed by 

To our knowledge, there is still lack of theoretical analysis  for the two dimensional quadratic termination of  the Dai--Yang method \cite{dai2006new2} whose stepsize is given by
\begin{equation} \label{Dai--Yang}
	\alpha_{k}^{P}=\frac{\left\|g_{k}\right\|_2}{\left\|A g_{k}\right\|_2}.
\end{equation}
A remarkable property of the Dai--Yang method is that $\alpha_{k}^{P}$ converges to the optimal stepsize $\frac{2}{\lambda_{1}+\lambda_{n}}$ (in the sense that it minimizes the modulus $\|I-\alpha A\|_{2}$, see \cite{dai2006new2,elman1994inexact}). Moreover, the Dai--Yang method is able to find the eigenvectors corresponding to $\lambda_{1}$ and $\lambda_{n}$.

%This stepsize converges to $\frac{2}{\lambda_{1}+\lambda_{n}}$, which  is in some sense an optimal stepsize because it minimizes the modulus $\|I-\alpha A\|_{2}$, see \cite{dai2006new2,elman1994inexact}.

In this paper, for a uniform analysis, we consider  equipping the family 
\begin{equation}\label{alphaAOfamily}
	\alpha_{k}=\frac{\left\|\psi(A)g_{k}\right\|_2}{\left\|\psi(A)A g_{k}\right\|_2}
\end{equation}
with the two dimensional quadratic termination property, which will be achieved by cooperating with
\begin{equation}\label{defalp}
	\tilde{\alpha}_{k}=\arg\max_{\alpha_{k} \in \mathbb{R}}~\alpha_{k+1}=\arg\max_{\alpha_{k} \in \mathbb{R}}~\frac{\left\|\psi(A)\left(I-\alpha_{k}A\right)g_{k}\right\|_2}{\left\|\psi(A)A \left(I-\alpha_{k}A\right)g_{k}\right\|_2}.
\end{equation}
The above strategy of maximizing the stepsize value in the next iteration has been employed in \cite{frassoldati2008new} for the SD method. However, the analysis in \cite{frassoldati2008new} can not be directly applied to the family \eqref{alphaAOfamily}. 
Clearly, $\alpha_{k}^{P}$ corresponds to the case $\psi(A) = I$ in \eqref{alphaAOfamily}.
We prove that each method in the family \eqref{alphaAOfamily} will asymptotically alternate in a two dimensional subspace associated with  the two eigenvectors corresponding to $\lambda_{1}$ and $\lambda_{n}$. In addition, for any given $\psi$, the stepsize \eqref{alphaAOfamily} tends to the above optimal stepsize as $k\rightarrow\infty$, and the eigenvectors corresponding to $\lambda_{1}$ and $\lambda_{n}$ can be obtained.  Then, we show that $\lim _{k \rightarrow \infty}\tilde{\alpha}_{k} =1 /\lambda _{n}$ for $n$ dimensional strictly convex quadratics. By adaptively taking $\alpha_{k}^{BB1}$ and  reusing  $\tilde{\alpha}_{k-1}$ for some iterations, we  propose a new method for quadratic minimization problems. It is proved that the proposed method is $R$-linearly convergent with the rate of $1-1/\kappa$, where $\kappa=\lambda_{n} / \lambda_{1}$ is the condition number of $A$. Our numerical comparisons with the BB1 \cite{barzilai1988two}, DY \cite{dai2005analysis}, SL (Alg.1  in \cite{sun2020new}), ABBmin2 \cite{frassoldati2008new}, SDC \cite{de2014efficient}, and MGC \cite{zou2021fast} methods for solving   unconstrained  random and non-random quadratic optimization  demonstrate that the proposed method is very efficient. Further, numerical experiments on quadratic problems whose Hessians are chosen from the SuiteSparse Matrix Collection \cite{davis2011university}  suggest that the proposed method  is very competitive with the above methods.

%Furthermore, we present a property of the stepsize which generalizes the Property A in \cite{dai2003alternate}, and establish a uniform $R$-linear convergence result for gradient methods satisfying the property. Consequently, the proposed method and many others are $R$-linearly convergent. Under proper conditions, we prove the $R$-factor can be bounded above by $1-1/\kappa$, where $\kappa=\lambda_{n} / \lambda_{1}$ is the condition number of $A$. 

%In Section \ref{sec:2}, we analyze the asymptotic behavior of the family \eqref{alphaAOfamily}. 

The paper is organized as follows. In Section \ref{sec:3}, we derive the new stepsize $\tilde{\alpha}_{k}$ and analyze its properties. The asymptotic behavior of the family \eqref{alphaAOfamily} is also analyzed. Our new algorithm for quadratic minimization problems as well as its $R$-linear convergence are presented in Section \ref{sec:4}. Section \ref{sec:5} presents some numerical comparisons of the proposed method and other successful gradient methods on solving  quadratic  problems. Finally, in Section \ref{sec:6} we give some concluding remarks.

\section{A new stepsize and its properties}
\label{sec:3}
%Text with citations \cite{RefB} and \cite{RefJ}.
%In this section,  we first analyze the asymptotic behavior of the family \eqref{alphaAOfamily} which is useful for analyzing the asymptotic spectral property  of the new stepsize $\tilde{\alpha}_{k}$. Then we derive the formula of $\tilde{\alpha}_{k}$ and present its  properties.
In this section, we derive the formula of $\tilde{\alpha}_{k}$ and analyze its  properties.

%\subsection{The stepsize}
To obtain $\tilde{\alpha}_{k}$, we consider
%By \eqref{alphaAOfamily},  $\alpha_{k+1}^{2}$  can be written  as 
\begin{equation}\label{delta13}	
%		\alpha_{k+1}^{2}= f_{1}(\alpha_{k})  f_{2}(\alpha_{k}) = 	\frac{\Delta_{1}(\alpha_{k}) }{\Delta_{3}(\alpha_{k})}:=F(\alpha_{k}).	
		F(\alpha_{k}):=\alpha_{k+1}^{2}=
		\frac{\left\|\psi(A)\left(I-\alpha_{k}A\right)g_{k}\right\|_2^{2}}{\left\|\psi(A)A \left(I-\alpha_{k}A\right)g_{k}\right\|_2^{2}}.
\end{equation}
The maximum value of $F(\alpha_{k})$ is achieved when $F'(\alpha_{k}) = 0$, which holds for any ${\alpha}_{k}$ satisfying
\begin{equation}\label{qpdF}
	\phi_{1} \alpha_{k}^{2}-\phi_{2} \alpha_{k}+\phi_{3}=0,
\end{equation}
where $\phi_{1}=c_{1} c_{4}-c_{2} c_{3}$, $\phi_{2}=c_{0} c_{4}-c_{2}^{2}$ and $\phi_{3}=c_{0} c_{3}-c_{1}c_{2}$
%\begin{equation*}\label{phi}
%	\phi_{1}=c_{1} c_{4}-c_{2} c_{3}, 	\qquad	\phi_{2}=c_{0} c_{4}-c_{2}^{2},	\qquad	\phi_{3}=c_{0} c_{3}-c_{1}c_{2},
%\end{equation*}
with
\begin{equation}\label{cj}
	c_{j} = g_{k}^{T} A^{j}  \psi^{2}(A)  g_{k},  \qquad j = 0,1,2,3,4.
\end{equation}

The following lemma guarantees that \eqref{qpdF} has two roots.
\begin{lemma}\label{Lemma1}
	Assume  $g_{k} \neq 0$ and $g_{k}$ is not parallel to $Ag_{k}$. Then, $\phi_{1}, \phi_{2}, \phi_{3} > 0$ and $\phi_{2}^2 - 4\phi_{1}\phi_{3}  > 0$.
%	\begin{equation}
%		\phi_{1}, \phi_{2}, \phi_{3} > 0,
%	\end{equation}
%	\begin{equation}\label{Lemma1.4}
%		\phi_{2}^2 - 4\phi_{1}\phi_{3}  > 0.
%	\end{equation}
\end{lemma}
\begin{proof}
It follows from the Cauchy-Schwartz inequality and \eqref{cj} that
\begin{equation*}	
		\begin{aligned}
			\left\|\psi(A)A^{\frac{j}{2}}g_{k}\right\|_2^2  \left\|\psi(A)A^{\frac{j+2}{2}}g_{k}\right\|_2^2
			&>\left(\left(\psi(A)A^{\frac{j}{2}}g_{k}\right)^T\left(\psi(A)A^{\frac{j+2}{2}}g_{k}\right)\right)^2\\	&=\left(g_{k}^{T}\psi^{2}(A)A^{j+1}g_{k}\right)^{2}
		\end{aligned}
\end{equation*}	
for $j\geq0$. That is, $c_{j} / c_{j+1} > c_{j+1}/ c_{j+2}$, which implies $\phi_{1}, \phi_{2}, \phi_{3} > 0.$

By direct calculation, we obtain $c_{1} \phi_{2} = c_{0}  \phi_{1} + c_{2} \phi_{3}$,
%	\begin{equation*}\label{Lemma1.2}
%		c_{1} \phi_{2} = c_{0}  \phi_{1} + c_{2} \phi_{3}.
%	\end{equation*}
which implies that
	\begin{equation*}\label{lemma1.5}
		c_{1}^{2}\phi_{2}^{2} =\left(c_{0}  \phi_{1} + c_{2} \phi_{3}\right)^{2}  \geq 4c_{0}c_{2}\phi_{1} \phi_{3}.
	\end{equation*}
Combining with  $c_{0}  c_{2}>c_{1}^{2}$, we have $\phi_{2}^2 - 4\phi_{1}\phi_{3}  > 0$. This completes the proof.\qed
\end{proof}

Using the square root law, we get the two roots of \eqref{qpdF} as 
\begin{equation}\label{alpha1}
	\tilde{\alpha}_{k} = \frac{\phi_{2} - \sqrt{\phi_{2}^{2}-4 \phi_{1} \phi_{3}}}{2 \phi_{1}}=\frac{2}{\frac{\phi_{2}}{\phi_{3}} + \sqrt{\left(\frac{\phi_{2}}{\phi_{3}}\right)^{2}-4 \frac{\phi_{1}}{\phi_{3}}}}
\end{equation}
and
\begin{equation}\label{alpha2}
	\hat{\alpha}_{k} = \frac{\phi_{2} + \sqrt{\phi_{2}^{2}-4 \phi_{1} \phi_{3}}}{2 \phi_{1}}=\frac{2}{\frac{\phi_{2}}{\phi_{3}} - \sqrt{\left(\frac{\phi_{2}}{\phi_{3}}\right)^{2}-4 \frac{\phi_{1}}{\phi_{3}}}}.
\end{equation}
It is easy to see that $\tilde{\alpha}_{k} = \arg\max_{\alpha_{k} \in \mathbb{R}}~\alpha_{k+1}$ and $\hat{\alpha}_{k} =\arg\min_{\alpha_{k} \in \mathbb{R}}~\alpha_{k+1}$.

Next theorem presents the two dimensional quadratic termination of the gradient method using $\alpha_{k}$ in \eqref{alphaAOfamily} and $\tilde{\alpha}_{k}$.
\begin{theorem}[Two dimensional quadratic termination]\label{Finite termination}
	Consider the gradient method \eqref{aa} for minimizing the two dimensional quadratic function \eqref{1a}. If the
	stepsize $\alpha_{k}$ is given by \eqref{alphaAOfamily} for all $k \neq k_{0}$ and   $ \alpha_{k_{0}}=\tilde{\alpha}_{k_{0}}$   at the $k_{0}$-th iteration where $k_{0} \geq 1$, it holds that $g_{k_{0}+i}=0$ for some $1 \leq i \leq 3$.
\end{theorem}
\begin{proof}
Without loss of generality, we assume that $A=\operatorname{diag}\{1, \lambda\}$ with $\lambda>0$. Let $g_{k}^{(1)}$ and $g_{k}^{(2)}$ be the first and second components of $g_{k}$, respectively.
	Notice that
\begin{equation*}
	c_{j}=g_{k}^{T} \psi^{2}(A) A^{j} g_{k}=\left(g_k^{(1)}\right)^{2} \psi^{2}(1)+\lambda^{j}\left(g_k^{(2)}\right)^{2} \psi^{2}(\lambda).
\end{equation*}
After direct calculation and simplification, we get
\begin{equation*}
	\begin{aligned}
		\phi_{1} &=\left(g_k^{(1)}\right)^{2}\left(g_k^{(2)}\right)^{2} \psi^{2}(1) \psi^{2}(\lambda)(\lambda-1)^{2}(\lambda+1) \lambda, \\
		\phi_{2} &=\left(g_k^{(1)}\right)^{2}\left(g_k^{(2)}\right)^{2} \psi^{2}(1) \psi^{2}(\lambda)(\lambda-1)^{2}(\lambda+1)^{2}, \\
		\phi_{3} &=\left(g_k^{(1)}\right)^{2}\left(g_k^{(2)}\right)^{2} \psi^{2}(1) \psi^{2}(\lambda)(\lambda-1)^{2}(\lambda+1).
	\end{aligned}
\end{equation*}	
	Therefore, we obtain
	\begin{equation*}
		\frac{\phi_{1}}{\phi_{3}} = \lambda \quad \text{and} \quad \frac{\phi_{2}}{\phi_{3}} = \lambda +1.
	\end{equation*}
Thus, from \eqref{alpha1} we know that $\tilde{\alpha}_{k} = 1/\lambda$ for all $k\geq 1$.  The conclusion follows immediately from $\alpha_{k_{0}}=\tilde{\alpha}_{k_{0}}$ and $g_{k+1}=\left(I-\alpha_{k}A\right)g_{k}$. We complete the proof.\qed
\end{proof}

%\subsection{Asymptotic behavior of the family \eqref{alphaAOfamily}}
In what follows, we shall prove that $\tilde{\alpha}_{k}$ converges to $1/\lambda_{n}$ under each method in the family \eqref{alphaAOfamily}. To this aim, we have to analyze the asymptotic behavior of the family \eqref{alphaAOfamily} first.

For convenience, we assume without loss of generality that the matrix $A$ is diagonal  with distinct eigenvalues, i.e.
\begin{equation}\label{diag}
	A=\mathrm{diag}\{\lambda_{1}, \lambda_{2}, \ldots, \lambda_{n}\}, \quad 0<\lambda_{1}<\lambda_{2}<\ldots<\lambda_{n}.
\end{equation}
Let $\left\{\xi_{1}, \xi_{2}, \ldots, \xi_{n}\right\}$ be the set of orthogonal eigenvectors associated with the eigenvalues $\left\{\lambda_{1}, \lambda_{2}, \ldots, \lambda_{n}\right\}$. 
Denoting by $\mu_{k}^{(i)}, i=1, \ldots, n$, the components of $g_{k}$ along $\xi_{i}$, i.e.
\begin{equation}\label{gk}
	g_{k}=\sum_{i=1}^{n} \mu_{k}^{(i)} \xi_{i}.
\end{equation}
It follows from \eqref{aa} and \eqref{gk} that
\begin{equation}\label{1c}
	g_{k+1}=\left(I-\alpha_{k}A\right)g_{k}=\prod_{j=0}^{k}\left(I-\alpha_{j} A\right) g_{0}=\sum_{i=1}^{n} \mu_{k+1}^{(i)} \xi_{i},
\end{equation}
%Then, substituting  \eqref{gk} into \eqref{1c} yields
where
\begin{equation}\label{mu}
	\mu_{k+1}^{(i)}=\left(1-\alpha_{k} \lambda_{i}\right) \mu_{k}^{(i)}=\mu_{0}^{(i)} \prod_{j=0}^{k}\left(1-\alpha_{j} \lambda_{i}\right).
\end{equation}

The following lemma is useful in our analysis.
%for analyzing the asymptotic behavior of the family \eqref{alphaAOfamily}.

\begin{lemma}\cite{dai2006new2}\label{lemmap*}	
	Let  $p$ be a vector in $\mathbb{R}^{n}$ such that (i)~$p^{(1)}>0$ and $p^{(n)}>0$  (ii) $p^{(1)} + p^{(n)} = 1$. Further assume that  $0<\lambda_{1}<\cdots<\lambda_{n}$. Consider
	a transformation $T$ such that
	\begin{equation*}\label{T}
		\left(Tp\right)^{(i)}=\frac{\left(\lambda_{i}-\gamma(p)\right)^{2} p^{(i)}}{\sum_{i}\left(\lambda_{i}-\gamma(p)\right)^{2} p^{(i)}}, \quad \text {where}\quad \gamma(p)=\sqrt{\sum_{i} \lambda_{i}^{2} p^{(i)}}.
	\end{equation*}
	Then
	\begin{equation*}
		\lim _{k \rightarrow \infty} T^{k}p= \begin{cases} h_{1}, & \text { if } i=1, \\ 0, & \text { if } i=2, \ldots, n-1, \\ h_{2}, & \text { if } i=n,\end{cases}
	\end{equation*}
	where 
	\begin{equation}\label{h1h2}
		h_{1}=\frac{\lambda_{1}+3 \lambda_{n}}{4\left(\lambda_{1}+\lambda_{n}\right)} \quad \text { and } \quad h_{2}=\frac{3 \lambda_{1}+ \lambda_{n}}{4\left(\lambda_{1}+\lambda_{n}\right)}.
	\end{equation}
\end{lemma}

Based on Lemma \ref{lemmap*}, we are able to show that each method in the family \eqref{alphaAOfamily} will zigzag in a two dimensional subspace which generalizes the results in \cite{dai2006new2}, where the case $\psi(A)=I$ (i.e. the Dai--Yang method) is considered.
\begin{theorem}\label{AOcon}
	Assume that the starting point $x_{0}$  is such that
	\begin{equation*}
		g_{0}^{T} \xi_{1} \neq 0 \quad \text { and }  \quad g_{0}^{T} \xi_{n} \neq 0.
	\end{equation*}
	Let $\{x_{k} \}$ be the sequence generated by applying a method in the family \eqref{alphaAOfamily}. Then
	\begin{equation}\label{alphakc1}
		\lim _{k \rightarrow \infty} \frac{\psi(\lambda_{i})\mu_{2k}^{(i)}}{\sqrt{{\sum_{j=1}^{n}\left(\psi(\lambda_{i})\mu_{2k}^{(i)}\right)^{2}}}}=\left\{
		\begin{array}{ll}
			\mathrm{sign}\left(\psi(\lambda_{1})\mu_{2k}^{(1)}\right) \sqrt{h_{1}}, & \text { if } i=1, \\
			0, & \text { if } i=2, \ldots, n-1, \\
			\mathrm{sign}\left(\psi(\lambda_{n})\mu_{2k}^{(n)}\right) \sqrt{h_{2}}, & \text { if } i=n,
		\end{array}\right.
	\end{equation}
	and
	\begin{equation}\label{alphakc2}
		\lim _{k \rightarrow \infty} \frac{\psi(\lambda_{i})\mu_{2k+1}^{(i)}}{\sqrt{\sum_{j=1}^{n}\left(\psi(\lambda_{i})\mu_{2k+1}^{(i)}\right)^{2}}}=
		\left\{\begin{array}{ll}
			\mathrm{sign}\left(\psi(\lambda_{1})\mu_{2k+1}^{(1)}\right) \sqrt{h_{1}}, & \text { if } i=1, \\
			0, & \text { if } i=2, \ldots, n-1, \\
			\mathrm{-sign}\left(\psi(\lambda_{n})\mu_{2k+1}^{(n)}\right) \sqrt{h_{2}}, & \text { if } i=n,
		\end{array}\right.
	\end{equation}
	where $h_{1}$ and $h_{2}$ are defined in \eqref{h1h2}.
	Further, the vectors\\
	$$
	\frac{\psi(A)g_{k}}{\| \psi(A)g_{k}\|_{2}}+\frac{\psi(A)g_{k+1}}{\left\|\psi(A)g_{k+1}\right\|_{2}}
	~~~\text{and}~~~
	\frac{\psi(A)g_{k}}{\| \psi(A)g_{k}\|_{2}}-\frac{\psi(A)g_{k+1}}{\left\|\psi(A)g_{k+1}\right\|_{2}}
	$$
	tend to be the eigenvectors corresponding to $\lambda_{1}$ and $\lambda_{n}$ of $A$, respectively.
\end{theorem}
\begin{proof}
	By  \eqref{gk}, we get
	\begin{equation}\label{psiAg}
		\psi(A)g_{k}=\sum_{i=1}^{n}\eta_{k}^{(i)} \xi_{i},
	\end{equation}
	where $\eta_{k}^{(i)} = \psi(\lambda_{i}) \mu_{k}^{(i)}$,
	%	\begin{equation*}\label{eta}
	%		\eta_{k}^{(i)} = \psi(\lambda_{i}) \mu_{k}^{(i)},
	%	\end{equation*}
	which together with \eqref{mu} yields
	\begin{equation}\label{eta+1}
		\eta_{k+1}^{(i)} = \left(1-\alpha_{k}\lambda_{i}\right)\eta_{k}^{i}.
	\end{equation}
	Defining the vector $p_{k}=\left(p_{k}^{(i)}\right)$ with
	\begin{equation}\label{pk}
		p_{k}^{(i)}=\frac{\left(\eta_{k}^{(i)}\right)^{2}}{\left\|\eta_{k}\right\|_{2}^{2}}
	\end{equation}
	and
	\begin{equation}\label{zeta}
		\gamma_{k}=	\alpha_{k}^{-1}=\frac{\left\|\psi(A)A g_{k}\right\|_{2}}{\left\|\psi(A)g_{k}\right\|_{2}}
		=\sqrt{\sum_{i} \lambda_{i}^{2} p_{k}^{(i)}}.
	\end{equation}
	We have from \eqref{eta+1}, \eqref{pk} and \eqref{zeta} that
	\begin{equation*}
		p_{k+1}^{(i)}=\frac{\left(\lambda_{i}-\gamma_{k}\right)^{2} p_{k}^{(i)}}{\sum_{i}\left(\lambda_{i}-\gamma_{k}\right)^{2} p_{k}^{(i)}}.
	\end{equation*}
	Clearly, according  to the definition of $p_{k}$, we get that $p_{k}^{(i)}\geq 0$ for all $i$ and
	\begin{equation*}
		\sum_{i} p_{k}^{(i)}=1, \quad \text { for all } k.
	\end{equation*}
	Let $p=p_{1} \in \mathbb{R}^{n}$, based on the above analysis and Lemma \ref{lemmap*}, we know that
	$
	\lim _{k \rightarrow \infty} p_{k}=\left(h_{1}, 0, \ldots, 0, h_{2}\right)^{T}
	$, where $h_{1}$ and $h_{2}$ are given in \eqref{h1h2}. 
	It follows from \eqref{eta+1} and $\lambda_{n}^{-1}<\alpha_{k}<\lambda_{1}^{-1}$ that
	\begin{equation}\label{sign}
		\mathrm{sign}\left(\psi(\lambda_{1})\mu_{k+1}^{(1)}\right)=\mathrm{sign}\left(\psi(\lambda_{1})\mu_{k}^{(1)}\right) ~ 
	\end{equation}
	and
	\begin{equation}\label{sign2}
		\mathrm{sign}\left(\psi(\lambda_{n})\mu_{k+1}^{(n)}\right)=-\mathrm{sign}\left(\psi(\lambda_{n})\mu_{k}^{(n)}\right).
	\end{equation}
	Thus, by Lemma \ref{lemmap*}, \eqref{sign} and \eqref{sign2}, we know that \eqref{alphakc1} and \eqref{alphakc2} hold. 
	
	Furthermore, combining \eqref{psiAg}, \eqref{sign} and \eqref{sign2}, we find that
	$$
	\lim _{k \rightarrow \infty} \frac{\psi(A)g_{k}}{\| \psi(A)g_{k}\|_{2}}+\frac{\psi(A)g_{k+1}}{\left\|\psi(A)g_{k+1}\right\|_{2}}=2 \mathrm{sign}\left(\psi(\lambda_{1})\mu_{2k}^{(1)}\right) \sqrt{h_{1}} \xi_{1}
	$$
	and
	$$
	\lim _{k \rightarrow \infty} \frac{\psi(A)g_{k}}{\| \psi(A)g_{k}\|_{2}}-\frac{\psi(A)g_{k+1}}{\left\|\psi(A)g_{k+1}\right\|_{2}}=\pm 2 \sqrt{h_{2}} \xi_{n}.
	$$
	This completes our proof.\qed
\end{proof}

From Theorem \ref{AOcon}, we have the following asymptotic result of the stepsize  \eqref{alphaAOfamily}.
\begin{corollary}\label{AOconalp}
	Under the conditions of Theorem \ref{AOcon}, for $\alpha_{k}$ in \eqref{alphaAOfamily} it holds that
	\begin{equation*}\label{alpconv}
		\lim _{k \rightarrow \infty} \alpha_{k}=\frac{2}{\lambda_1+\lambda_n}.
	\end{equation*}
\end{corollary}

Next theorem shows that $\tilde{\alpha}_{k}$ converges to $1/\lambda_{n}$ under each method in the family \eqref{alphaAOfamily}.
\begin{theorem}\label{theorem5}	
	Under the conditions of Theorem \ref{AOcon}, let $\{g_{k}\}$ be the sequence generated by applying a method in the family \eqref{alphaAOfamily} to minimize the $n$-dimensional quadratic function \eqref{1a}. Then $\lim _{k \rightarrow \infty}\tilde{\alpha}_{k}=1/\lambda_{n}$.	
\end{theorem}
\begin{proof}
From \eqref{1c} and \eqref{cj}, we obtain
	\begin{equation}\label{cj/mu}
		c_{j} = g_{k}^{T} A^{j}\psi^{2}(A)  g_{k}
		=  \sum_{i=1}^{n} \lambda_{i}^{j}\left(\psi(\lambda_{i})\mu_{k}^{(i)}\right)^{2}.
	\end{equation}
	When $k$ is odd,  by the definition of $\phi_{1}$, \eqref{alphakc2} and \eqref{cj/mu}, we get	
	\begin{align}\label{phi1odd}
		\lim _{k \rightarrow \infty} \frac{\phi_{1}}{c_{0}^{2}} \nonumber&= \lim _{k \rightarrow \infty} \left(\frac{c_{1}}{c_{0}} \frac{c_{4}}{c_{0}}-\frac{c_{2}}{c_{0}} \frac{c_{3}}{c_{0}}\right)\\
		\nonumber	& = \lim _{k \rightarrow \infty}\left[\sum_{i=1}^{n} \lambda_{i} \frac{\left(\psi(\lambda_{i})\mu_{k}^{(i)}\right)^{2}}{\sum_{s=1}^{n}\left(\psi(\lambda_{s})\mu_{k}^{(s)}\right)^{2}} \cdot \sum_{i=1}^{n} \lambda_{i}^{4} \frac{\left(\psi(\lambda_{i})\mu_{k}^{(i)}\right)^{2}}{\sum_{s=1}^{n}\left(\psi(\lambda_{s})\mu_{k}^{(s)}\right)^{2}} \right] \\
		\nonumber	&\quad -\lim _{k \rightarrow \infty}\left[\sum_{i=1}^{n} \lambda_{i}^2 \frac{\left(\psi(\lambda_{i})\mu_{k}^{(i)}\right)^{2}}{\sum_{s=1}^{n}\left(\psi(\lambda_{s})\mu_{k}^{(s)}\right)^{2}} \cdot \sum_{i=1}^{n} \lambda_{i}^{3} \frac{\left(\psi(\lambda_{i})\mu_{k}^{(i)}\right)^{2}}{\sum_{s=1}^{n}\left(\psi(\lambda_{s})\mu_{k}^{(s)}\right)^{2}} \right] \\
		\nonumber	& =\left(h_{1}\lambda_{1}+  h_{2}\lambda_{n}\right)\left(h_{1}\lambda_{1}^{4}+  h_{2}\lambda_{n}^{4}\right) -\left(h_{1}\lambda_{1}^{2}+  h_{2}\lambda_{n}^{2}\right)
		\left(h_{1}\lambda_{1}^{3}+  h_{2}\lambda_{n}^{3}\right)\\
		& =h_{1}h_{2}\lambda_{1}\lambda_{n}\left(\lambda_{n}-\lambda_{1}\right)^{2}\left(\lambda_{1}+\lambda_{n}\right).
	\end{align}
	Similarly, we have 
	\begin{equation}\label{phi2odd}
		\lim _{k \rightarrow \infty}  \frac{\phi_{2}}{c_{0}^{2}} =  h_{1}h_{2}\left(\lambda_{n}-\lambda_{1}\right)^{2}\left(\lambda_{1}+\lambda_{n}\right)^{2} 
	\end{equation}
	and
	\begin{equation}\label{phi3odd}
		\lim _{k \rightarrow \infty}  \frac{\phi_{3}}{c_{0}^{2}} =    h_{1}h_{2}\left(\lambda_{n}-\lambda_{1}\right)^{2}\left(\lambda_{1}+\lambda_{n}\right).
\end{equation}
	Combining \eqref{phi1odd},  \eqref{phi2odd} and \eqref{phi3odd}, we obtain
	\begin{equation} \label{phi1phi2phi3}
		\lim _{k \rightarrow \infty} \frac{\phi_{1}}{\phi_{3}} = \lambda_{1}\lambda_{n} \quad  \text{ and } \quad	\lim _{k \rightarrow \infty} \frac{\phi_{2}}{\phi_{3}} = \lambda_{1}+\lambda_{n}.
	\end{equation}
	It follows from \eqref{phi1phi2phi3} and \eqref{alpha1} that
	$\lim _{k \rightarrow \infty} \tilde{\alpha}_{k}=1/\lambda_{n}$. When  $k$ is even, we get the desired result in the same manner as above.  This completes our proof.
	\qed
\end{proof}

\begin{figure}[htp!b]\label{Figure 1}
	\centering
	\		\includegraphics[width=3in]{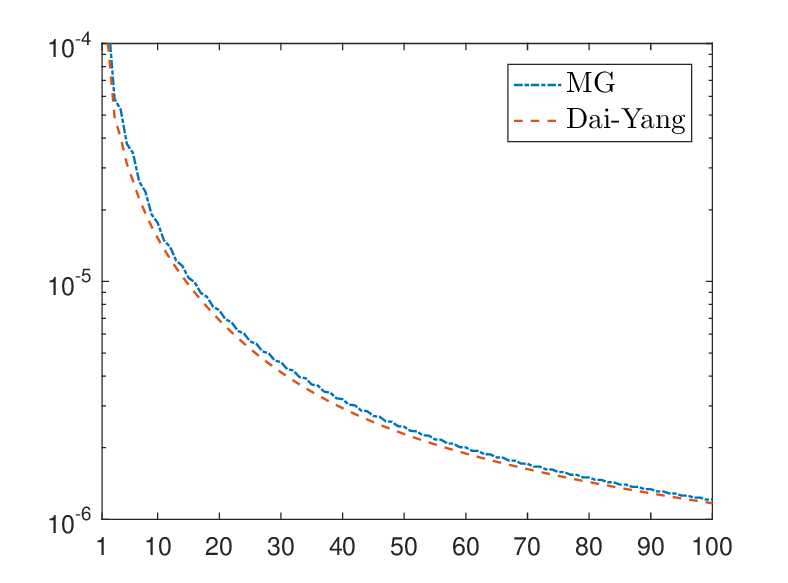}	
	\caption{ Problem \eqref{problem1} with $n = 1000$:  convergence history of the sequence $\{|\tilde{\alpha}_{k}-1/\lambda_{n}|\}$ for
		the first 100 iterations of the MG and Dai--Yang methods.}
\end{figure}

From \eqref{phi1phi2phi3}, we see that  $\phi_{1}/\phi_{3}$ and  $\phi_{2}/\phi_{3}$ are independent of $\psi(A)$. The following example shows that  $\tilde{\alpha}_{k}$  converges to $1 / \lambda_{n}$ under both the  MG  and Dai--Yang methods. In particular, we applied the  MG  and Dai--Yang methods to the quadratic function \eqref{1a} with
\begin{equation}\label{problem1}
	A=\mathrm{diag}\left\{a_{1}, a_{2}, \ldots, a_{n}\right\} \quad \text { and } \quad b=0,
\end{equation}
where $a_{1}=1,~a_{n}=n$ and $a_{i}$ was randomly generated in $(1, n)$, $i=2, \ldots, n-1$. The starting point was set to the vector of all ones.  From Figure 1,  we see that  $\tilde{\alpha}_{k}$   approximates  $1 / \lambda_{n}$  with satisfactory accuracy after a small number of iterations under both  the  MG  and Dai--Yang methods.

\section{A new gradient method}\label{sec:4}
In this section, we propose a new algorithm for unconstrained quadratic optimization and present its $R$-linear convergence result.

%\subsection{The algorithm}
Extensive studies  point out that  adaptively choosing a short stepsize or $\alpha_{k}^{BB1}$ at each iteration is numerically better than the original BB method, see for example \cite{crisci2020spectral,dai2005projected,dai2006cyclic,frassoldati2008new,huang2020acceleration,zhou2006gradient}. Now we show that the new stepsize $\tilde{\alpha}_{k}$ is a short one.
\begin{theorem}\label{theoremalpha}
	Under the conditions of Lemma \ref{Lemma1}, it holds that
	$\tilde{\alpha}_{k} < c_{2}/c_{3}$.
\end{theorem}
\begin{proof}
	According to  \eqref{cj}  and \eqref{alpha1}, we  get
	\begin{equation}\label{theorem1.1}
		\frac{c_{2}}{c_{3}} - \tilde{\alpha}_{k} = \frac{c_{2} \sqrt{\phi_{2}^{2}-4\phi_{1}\phi_{3}} + c_{2}\phi_{2} - 2 c_{3} \phi_{3}}{c_{3}\left(\phi_{2} + \sqrt{\phi_{2}^{2}-4\phi_{1}\phi_{3}}\right)}.
	\end{equation}
	If $c_{2}\phi_{2} - 2 c_{3} \phi_{3} > 0$, we have  $c_{2}/c_{3} - \tilde{\alpha}_{k} > 0$; otherwise, it follows from \eqref{theorem1.1} that
	%	By the definitions of $\phi_{1},\phi_{2} \text{ and } \phi_{3}$, we get
	%	\begin{equation}\label{Lemma1.3}
	%		c_{3}\phi_{2} = c_{2}\phi_{1} +c_{4} \phi_{3}.
	%	\end{equation}
	%	Combining \eqref{theorem1.1}, \eqref{theorem1.2} and  \eqref{Lemma1.3}, we have
	\begin{equation*}
		\begin{aligned}
			\frac{c_{2}}{c_{3}} - \tilde{\alpha}_{k} &=\frac{4\left(c_{2}c_{3}\phi_{2}\phi_{3}-(c_{3}\phi_{3})^2-c_{2}^2\phi_{1}\phi_{3}\right)}{c_{3}\Delta\left(\phi_{2} + \sqrt{\phi_{2}^{2}-4\phi_{1}\phi_{3}}\right)}\\
			&=\frac{4\phi_{3}^{2}\left(c_{2}c_{4}-c_{3}^2\right)}{c_{3}\Delta\left(\phi_{2}+\sqrt{\phi_{2}^{2}-4\phi_{1}\phi_{3}}\right)},
		\end{aligned}
	\end{equation*}
%	\begin{equation*}\label{theorem1.2}
%	\Delta = c_{2} \sqrt{\phi_{2}^{2}-4\phi_{1}\phi_{3}} -\left( c_{2}\phi_{2} - 2 c_{3}\phi_{3}\right)>0,
%\end{equation*}
	where the last equality is due to $c_{3}\phi_{2} = c_{2}\phi_{1} +c_{4} \phi_{3}$ and $\Delta = c_{2} \sqrt{\phi_{2}^{2}-4\phi_{1}\phi_{3}} -\left( c_{2}\phi_{2} - 2 c_{3}\phi_{3}\right)>0$.
	From $\phi_{1}, \phi_{3}>0$ and $c_{2}c_{4}-c_{3}^2 > 0$, we know that $\tilde{\alpha}_{k}<c_{2}/c_{3}$ holds. This completes our proof.\qed
\end{proof}

%From Lemma \ref{Lemma1} and Theorem \ref{theoremalpha}, we know that $\tilde{\alpha}_{k}\leq \alpha_{k}^{MG}$  when $\psi(A)=I$.

Based on the above analysis, we can develop gradient method using $\tilde{\alpha}_{k}$ and $\alpha_{k}^{BB1}$ in an adaptive way. Notice that reusing the retard short stepsize for some iterations could reduce the computational cost and yield better performance, see \cite{de2014efficient,huang2020gradient,sun2020new,yuan2008step,zou2021fast}. So, we suggest to combine the adaptive and cyclic schemes with  $\alpha_{k}^{BB1}$ and  $\tilde{\alpha}_{k-1}$. In particular, our method reuses $\tilde{\alpha}_{k-1}$ for $r$ iterations when $\alpha_{k}^{\mathrm{BB} 2} / \alpha_{k}^{\mathrm{BB} 1}<\tau$ for some $\tau \in (0,1)$; otherwise, we set $\alpha_{k} = \alpha_{k}^{BB1}$. We use $t$ as the index to keep track of the number of short stepsizes  chosen during the iterative process. The stepsize for our algorithm is summarized as
\begin{equation}\label{alphagk}
	\alpha_{k}= \begin{cases} \tilde{\alpha}_{k-1}, & \text { if } \bmod (t, r)=0~\text{and} ~\alpha_{k}^{\mathrm{BB} 2} / \alpha_{k}^{\mathrm{BB} 1}<\tau; \\ \alpha_{k}^{\mathrm{BB} 1}, & \text { if } \bmod (t, r)=0~\text{and} ~\alpha_{k}^{\mathrm{BB} 2} / \alpha_{k}^{\mathrm{BB} 1} \geq\tau; \\ \alpha_{k-1}, & \text { otherwise.}
	\end{cases}
\end{equation}
We mention that the parameter $\tau$ can be chosen dynamically as \cite{huang2020equipping}. However, in our test the dynamic scheme does not show much evidence over the above fixed one.  
Our method  is formally presented in Algorithm \ref{Al1}. 
\begin{algorithm}[htp!b] 
	\caption{A gradient method for unconstrained quadratic optimization} 
	\label{Al1}
%	\begin{algorithmic}[1] 
	Choose $x_{0}\in\mathbb{R}^{n}$, $\epsilon,~\tau\in(0,1),~r \in \mathbb{N}$. Set $\alpha_{0}=\alpha_{0}^{SD},~k:=0,~t:=0.$\\	 
		\While{ $\left\|g_{k}\right\|>\epsilon$}{
		$x_{k+1}=x_{k}-\alpha_{k} g_{k}$\\	
		\eIf{\rm{mod}$(t,r) = 0$}
		{
			\eIf{$\alpha_{k+1}^{\mathrm{BB} 2} / \alpha_{k+1}^{\mathrm{BB} 1}<\tau$}{
				compute $\tilde{\alpha}_{k}$\\
				$\alpha_{k+1} = \tilde{\alpha}_{k}$\\
				$t:= t+1$
			}{
			$\alpha_{k+1} = \alpha_{k+1}^{\mathrm{BB} 1}$\\
		}
	}{
			$\alpha_{k+1} = \alpha_{k}$\\
			$t:= t +1$
		}
		$k:=k+1$
	}
%	\end{algorithmic} 
\end{algorithm}

To establish $R$-linear convergence of  Algorithm \ref{Al1}, we first show the boundedness of $\tilde{\alpha}_{k}$.

Let us denote 
\begin{equation}\label{Falpha1}
	f_{1}(\alpha)= \frac{g_{k}^{T} \left(I-\alpha A\right)^{T}\psi^{2}(A)\left(I-\alpha A\right) g_{k}}{g_{k}^{T}\left(I-\alpha A\right)^{T} \psi^{2}(A) A \left(I-\alpha A\right)g_{k}}
	%	= \frac{\Delta_{1}(\alpha)}{\Delta_{2}(\alpha)}
\end{equation}
and
\begin{equation}\label{Falpha2}
	f_{2}(\alpha) = \frac{g_{k}^{T} \left(I-\alpha A\right)^{T}\psi^{2}(A) A \left(I-\alpha A\right) g_{k}}{g_{k}^{T}\left(I-\alpha A\right)^{T} \psi^{2}(A) A^{2} \left(I-\alpha A\right)g_{k}}.
	%	= \frac{\Delta_{2}(\alpha) }{\Delta_{3}(\alpha)},
\end{equation}
Using the same arguments as \eqref{qpdF} and Lemma \ref{Lemma1}, we know that both $f_{1}'(\alpha)=0$ and $f_{2}'(\alpha)=0$ have two roots. In addition, the roots of $f_{1}'(\alpha)=0$, say $\tilde{\beta}_{k}$ and $\hat{\beta}_{k}$ with $\tilde{\beta}_{k}\leq\hat{\beta}_{k}$, are the solutions of
\begin{equation}\label{betak}
	\phi_{4}\beta^2-\phi_{3}\beta+\phi_{5}=0,
\end{equation}
and the roots of $f_{2}'(\alpha)=0$, say $\tilde{\gamma}_{k}$ and $\hat{\gamma}_{k}$ with $\tilde{\gamma}_{k}\leq\hat{\gamma}_{k}$, satisfy
\begin{equation}\label{gamak}
	\phi_{6}\gamma^2-\phi_{1}\gamma+\phi_{4}=0,
\end{equation}
where $\phi_{1}$, $\phi_{3}$ are defined in the former section, and $\phi_{4}=c_{1} c_{3}-c_{2}^{2}$, $\phi_{5}=c_{0} c_{2}-c_{1}^{2}$ and $\phi_{6}=c_{2} c_{4}-c_{3}^{2}$. Recall that $F(\alpha)=f_{1}(\alpha)  f_{2}(\alpha)$.  Since $f_{1}(\alpha)>0$ and $f_{2}(\alpha)>0$, any root of $F'(\alpha) = f_{1}'(\alpha)f_{2}(\alpha) + f_{1}(\alpha)f_{2}'(\alpha)=0$ yields
\begin{equation*}
	f_{1}'(\alpha)=f_{2}'(\alpha)=0\quad \text{or} \quad f_{1}'(\alpha) f_{2}'(\alpha)<0.
\end{equation*}
If the root is such that $f_{1}'(\alpha)\leq 0$ and $f_{2}'(\alpha)\geq 0$, we must have
\begin{equation}\label{bdalp1}
	\alpha\in[\tilde{\beta}_{k},\hat{\beta}_{k}]
	\quad \text{and} \quad
	\alpha\in(-\infty,\tilde{\gamma}_{k}]\cup [\hat{\gamma}_{k},+\infty).
\end{equation}
%where $\tilde{\beta}_{k}$ and $\hat{\beta}_{k}$ with $\tilde{\beta}_{k}\leq\hat{\beta}_{k}$ are roots of $f_{1}'(\alpha)=0$, and $\tilde{\gamma}_{k}$ and $\hat{\gamma}_{k}$ with $\tilde{\gamma}_{k}\leq\hat{\gamma}_{k}$ are roots of $f_{2}'(\alpha)=0$. 
Similarly, for the case $f_{1}'(\alpha)\geq 0$ and $f_{2}'(\alpha)\leq 0$, we get
\begin{equation}\label{bdalp2}
	\alpha\in[\tilde{\gamma}_{k},\hat{\gamma}_{k}]
	\quad \text{and} \quad
	\alpha\in(-\infty,\tilde{\beta}_{k}]\cup [\hat{\beta}_{k},+\infty).
\end{equation}

Now we are ready to prove that $\tilde{\alpha}_{k}$ is bounded between $1/\lambda_{n}$ and $1/\lambda_{1}$.
\begin{theorem}\label{alphalambda}
	Under the conditions of Lemma \ref{Lemma1}, it holds that $1/\lambda_{n}\leq\tilde{\alpha}_{k} \leq 1/\lambda_{1}$.
	%	$$\frac{1}{\lambda_{n}}\leq\tilde{\alpha}_{k}\leq \frac{1}{\lambda_{1}}.$$
\end{theorem}
\begin{proof}
	Since $F'(\tilde{\alpha}_{k}) =0$, when $f_{1}'(\tilde{\alpha}_{k})\leq0$, it follows from \eqref{bdalp1} that
	%(a) if  $f_{1}'(\alpha)\leq0$, then
	\begin{equation*}
		\tilde{\beta}_{k}\leq\tilde{\alpha}_{k} \leq \min\{\tilde{\gamma}_{k},\hat{\beta}_{k}\}
		\quad\text{or}\quad 	\max\{\hat{\gamma}_{k},\tilde{\beta}_{k}\} \leq \tilde{\alpha}_{k} \leq \hat{\beta}_{k}.
	\end{equation*}
	For the case  $f_{1}'(\tilde{\alpha}_{k})\geq 0$, by \eqref{bdalp2}, we get
	%(b) if  $f_{1}'(\alpha)\geq 0$, then
	\begin{equation*}
		\tilde{\gamma}_{k} \leq \tilde{\alpha}_{k} \leq \min\{\tilde{\beta}_{k},\hat{\gamma}_{k}\}
		\quad\text{or}\quad 	\max\{\hat{\beta}_{k},\tilde{\gamma}_{k}\} \leq \tilde{\alpha}_{k} \leq \hat{\gamma}_{k}.
	\end{equation*}
	Thus, we only need to prove that $1/\lambda_{n}\leq	\hat{\beta}_{k},~	\tilde{\beta}_{k},~\hat{\gamma}_{k},~	\tilde{\gamma}_{k} \leq 1/\lambda_{1}$.
	%	$$ \frac{1}{\lambda_{n}}\leq	\hat{\beta}_{k},~	\tilde{\beta}_{k},~\hat{\gamma}_{k},~	\tilde{\gamma}_{k} \leq \frac{1}{\lambda_{1}}. $$

	To proceed, we first show that 
	\begin{equation}\label{tildehat}
		\hat{\beta}_{k}=f_{1}(\tilde{\beta}_{k}).
	\end{equation}
	Since $f_{1}'(\tilde{\beta}_{k}) = 0$, we have
	\begin{equation*}
		(c_{1}-2  c_{2}\tilde{\beta}_{k}+ c_{3}\tilde{\beta}_{k}^{2})(-c_{1}+ c_{2}\tilde{\beta}_{k})-(c_{0}-2  c_{1}\tilde{\beta}_{k}+ c_{2}\tilde{\beta}_{k}^{2})(-c_{2}+ c_{3}\tilde{\beta}_{k})=0,
	\end{equation*}
	which implies that
	\begin{align*}
		f_{1}(\tilde{\beta}_{k})&=\frac{c_{0}-2  c_{1}\tilde{\beta}_{k}+c_{2}\tilde{\beta}_{k}^{2} }{c_{1}-2 c_{2} \tilde{\beta}_{k}+ c_{3}\tilde{\beta}_{k}^{2}}=\frac{c_{2} \tilde{\beta}_{k}-c_{1}}{c_{3} \tilde{\beta}_{k}-c_{2}}\\
		&=\hat{\beta}_{k} \frac{c_{2} \tilde{\beta}_{k}-c_{1}}{c_{3} \hat{\beta}_{k} \tilde{\beta}_{k}-c_{2} \hat{\beta}_{k}}=\hat{\beta}_{k} \frac{c_{2} \tilde{\beta}_{k}-c_{1}}{c_{3} \frac{\phi_{5}}{\phi_{4}}-c_{2} \hat{\beta}_{k}},
	\end{align*}
	where the last equality comes from the fact that $\tilde{\beta}_{k}\hat{\beta}_{k}=\phi_{5}/\phi_{4}$.
	%It follows from \eqref{betak} that $\tilde{\beta}_{k}\hat{\beta}_{k}=\phi_{5}/\phi_{4}$.
	%Hence
	%\begin{equation*}
	%	f_{1}(\tilde{\beta}_{k})=\hat{\beta}_{k} \frac{c_{2} \tilde{\beta}_{k}-c_{1}}{c_{3} \hat{\beta}_{k} \tilde{\beta}_{k}-c_{2} \hat{\beta}_{k}}=\hat{\beta}_{k} \frac{c_{2} \tilde{\beta}_{k}-c_{1}}{c_{3} \frac{\phi_{5}}{\phi_{4}}-c_{2} \hat{\beta}_{k}}.
	%\end{equation*}
	In order to prove \eqref{tildehat}, we are suffice to show that 
	$$c_{2} \tilde{\beta}_{k}-c_{1} = c_{3} \frac{\phi_{5}}{\phi_{4}}-c_{2} \hat{\beta}_{k},$$
	that is, 
	$$c_{2}\phi_{4} (\tilde{\beta}_{k}+\hat{\beta}_{k}) =c_{1}\phi_{4}+ c_{3}\phi_{5},$$ 
	which holds due to $\tilde{\beta}_{k}+\hat{\beta}_{k}=\phi_{3}/\phi_{4}$ and the definitions of $c_{1}, c_{2}, c_{3}, \phi_{3}, \phi_{4}, \phi_{5}$.
	%\begin{equation*}
	%	c_{2} \tilde{\beta}_{k}-c_{1} = c_{3} \frac{\phi_{5}}{\phi_{4}}-c_{2} \hat{\beta}_{k}.
	%\end{equation*}
	%Combining   \eqref{betak}  and $c_{2}\phi_{3}=c_{3}\phi_{5}+c_{1}\phi_{4}$  yields
	%\begin{equation*}
	%	\begin{aligned}
	%		c_{2} \tilde{\beta}_{k}-c_{1} &=c_{2} \frac{\phi_{3}-\sqrt{\phi_{3}^{2}-4\phi_{4}\phi_{5}}}{2 \phi_{4}}-c_{1} \\
	%		&=\frac{c_{2} \phi_{3}-c_{2} \sqrt{\phi_{3}^{2}-4\phi_{4}\phi_{5}}-2 c_{1}\phi_{4} }{2 \phi_{4} } \\
	%		&=\frac{c_{3} \phi_{5}-c_{1} \phi_{4}-c_{2} \sqrt{\phi_{3}^{2}-4 \phi_{4}\phi_{5}}}{2 \phi_{4}}.
	%	\end{aligned}
	%\end{equation*}
	%Similarly, we obtain
	%$$
	%c_{3} \frac{\phi_{5}}{\phi_{4}}-c_{2} \hat{\beta}_{k}=\frac{c_{3} \phi_{5}-c_{1} \phi_{4}-c_{2} \sqrt{\phi_{3}^{2}-4 \phi_{4}\phi_{5}}}{2 \phi_{4}}=c_{2} \tilde{\beta}_{k}-c_{1}.
	%$$
	Thus, \eqref{tildehat} is true. It follows from the definition of $f_1$ and the Rayleigh's quotient property that $1/\lambda_{n} \leq\hat{\beta}_{k}=f_{1}(\tilde{\beta}_{k})\leq 1/\lambda_{1}$.
	%$$\frac{1}{\lambda_{n}} \leq\hat{\beta}_{k}=f_{1}(\tilde{\beta}_{k})\leq \frac{1}{\lambda_{1}}.$$
	Using the same arguments as above, we get
	$1/\lambda_{n} \leq\tilde{\beta}_{k},~\tilde{\gamma}_{k},~\hat{\gamma}_{k} \leq 1/\lambda_{1}$. This completes the proof.
	\qed
\end{proof}

In \cite{huangdl2022r}, the authors prove that any gradient method with stepsizes satisfying the following Property B has $R$-linear convergence rate $1-\lambda_{1}/M_1$ which implies a $1-1/\kappa$ rate when $M_1\leq\lambda_n$. Similar results for gradient methods satisfying the Property A in \cite{dai2003alternate} can be found in \cite{huang2021On}. However, a stepsize satisfies Property B may not meets the conditions of Property A.

\noindent\textbf{Property~B:} We
say that the stepsize $\alpha_{k}$ has Property B if there exist an integer $m$ and positive constant $M_{1} \geq \lambda_{1}$ such that\\
(i) $\lambda_{1} \leq \alpha_{k}^{-1} \leq M_{1}$;\\
(ii) for some real analytic function $\psi$ defined as \eqref{psi} and $v(k) \in\{k,k-1,\ldots,\max\{k-m+1,0\}\}$,
\begin{equation}\label{alphack} 
\alpha_{k} \leq \frac{g_{v(k)}^{T} \psi^{2}(A) g_{v(k)}}{g_{v(k)}^{T} A\psi^{2}(A) g_{v(k)}}.
\end{equation}

Clearly, $\alpha_k^{BB1}$ satisfies Property B with $M_1=\lambda_n$, $m=2$ and $\psi^{2}(A)=I$, and the new stepsize $\tilde{\alpha}_{k}$ satisfies Property B with $M_1=\lambda_n$. So, we are able to show $R$-linear convergence of  Algorithm \ref{Al1} with $1-1/\kappa$ rate as Theorem 2 in \cite{huangdl2022r}. For completeness, we include the proof here.

\begin{theorem}\label{converth2}
Suppose that the sequence $\left\{g_{k}\right\}$ is generated by Algorithm \ref{Al1}. Then either $g_{k} = 0$ for some finite $k$ or the sequence  $\left\{\left\|g_{k}\right\|_{2}\right\}$ converges to zero $R$-linearly in the sense that
%	\begin{equation}\label{Citheta3}
%		|g_{k}^{(i)}| \leq \widetilde{C}_{i}\tilde{\theta}^{k}, \qquad i=1,2, \ldots, n,
%	\end{equation}
%	where  $\tilde{\theta}=1-1/\kappa$ and
%	\begin{equation*}
%	\left\{\begin{array}{l}
%		\widetilde{C}_{1}=|g_{0}^{(1)}|; \\
%		\widetilde{C}_{i}=\displaystyle\max\left\{|g_{0}^{(i)}|,\frac{|g_{1}^{(i)}|}{\tilde{\theta}}, \ldots,\frac{|g_{m-1}^{(i)}|}{\tilde{\theta}^{m-1}},\frac{\max \{\tilde{\sigma}_{i}, \tilde{\sigma}_{i}^{m}\}}{\tilde{\theta}^{m}\psi(\lambda_{i})}\sqrt{\sum_{j=1}^{i-1} \psi(\lambda_{j})^{2}\widetilde{C}_{j}^{2}}\right\}
%	\end{array}\right.
%\end{equation*}
%with $\tilde{\sigma}_{i}=
%\max \left\{\frac{\lambda_{i}}{\lambda_{1}}-1,1-\frac{\lambda_{i}}{\lambda_{n}}\right\}$.
	\begin{equation}\label{Citheta}
		|g_{k}^{(i)}| \leq C_{i}\theta^{k}, \qquad i=1,2, \ldots, n,
	\end{equation}
	where  ${\theta}=1-1/\kappa$ and
\begin{equation*}
	\left\{\begin{array}{l}
	C_{1}=|g_{0}^{(1)}|; \\
	C_{i}=\displaystyle\max\left\{|g_{0}^{(i)}|,\frac{|g_{1}^{(i)}|}{\theta}, \ldots,\frac{|g_{r-1}^{(i)}|}{\theta^{r-1}},\frac{\max \{\sigma_{i}, \sigma_{i}^{r}\}}{\theta^{r}\psi(\lambda_{i})}\sqrt{\sum_{j=1}^{i-1}\psi^{2}(\lambda_{j}) C_{j}^{2}}\right\},\\
	\quad i = 2,3, \ldots n,
	\end{array}\right.
\end{equation*}
with $\sigma_{i}=\max \left\{\frac{\lambda_{i}}{\lambda_{1}}-1,1-\frac{\lambda_{i}}{\lambda_{n}}\right\}$.
\end{theorem}

\begin{proof}
%	It follows from Theorem \ref{alphalambda} and the Cauchy-Schwartz inequality that $1 / \lambda_{1} \leq \tilde{\alpha}_{k-1}^{-1}, (\alpha_{k}^{BB1})^{-1} \leq$ $1 / \lambda_{n}$. Thus, (i) of Property B holds with $M = \lambda_{n}$. 	
It follows from $g_{k+1}^{(i)}=\left(1-\alpha_{k} \lambda_{i}\right) g_{k}^{(i)}$ and (i) of Property B that
\begin{equation}\label{gksigma}
\left|g_{k+1}^{(i)}\right| \leq \sigma_{i}\left|g_{k}^{(i)}\right|, \quad i=1,2, \ldots, n.
\end{equation}
%For $i=1$,  we obtain
%\begin{equation*}
%\left|g_{k}^{(1)}\right| \leq \theta\left|g_{k-1}^{(1)}\right| \leq C_{1} \theta^{k}.
%\end{equation*}
Clearly, \eqref{Citheta} holds for $i=1$. 

In what follows, we prove \eqref{Citheta} by induction on $i$. Assume that \eqref{Citheta} holds for all $1 \leq i \leq L-1$ with $L \in\{2, \ldots, n\}$. When $i=L$, it follows from the definition of $C_{L}$  that \eqref{Citheta} holds for $k = 0, 1,\ldots, r-1$. Assume by contradiction that \eqref{Citheta} does not hold for $k\geq r$ when $i=L$. Let $\hat{k}\geq r$ be the minimal index such that $\left|g_{\hat{k}}^{(L)}\right| > C_{L} \theta^{\hat{k}}$. Then, if $\alpha_{\hat{k}-1} \lambda_{L} \leq 1$, we get
\begin{equation*}
 \left|g_{\hat{k}}^{(L)}\right|=\left(1-\alpha_{\hat{k}-1} \lambda_{L}\right)\left|g_{\hat{k}-1}^{(L)}\right| \leq \theta\left|g_{\hat{k}-1}^{(L)}\right| \leq C_{L} \theta^{\hat{k}},
\end{equation*}
 which contradicts our assumption. Thus we must have $\alpha_{\hat{k}-1} \lambda_{L} > 1$, which combines with \eqref{gksigma} and $\theta<1$ gives that
 \begin{align*}\label{ineqGeps}
 \psi(\lambda_{L})|g_{\hat{k}-j}^{(L)}| 
 &\geq \frac{\psi(\lambda_{L})|g_{\hat{k}}^{(L)}|}{{\sigma}_{L}^{j}}>\frac{\psi(\lambda_{L})C_{L} \theta^{\hat{k}}}{{\sigma}_{L}^{j}} 
 \geq \frac{\psi(\lambda_{L})C_{L} \theta^{\hat{k}}}{\max \{\sigma_{L}, \sigma_{L}^{r}\}} 
 \\
 & \geq \theta^{\hat{k}-r}\sqrt{\sum_{i=1}^{L-1}\psi^{2}(\lambda_{i}) C_{i}^{2}} \nonumber
 \geq \theta^{\hat{k}-r}\sqrt{ \sum_{i=1}^{L-1}\left(\psi(\lambda_{i})g_{\hat{k}-j}^{(i)}\right)^{2} \theta^{2(j-\hat{k})}}  \nonumber\\
 &=\theta^{j-r} \sqrt{G(\hat{k}-j, L-1)} 
 \geq \sqrt{G(\hat{k}-j, L-1)},\quad j \in[1, r],
 \end{align*}
%\begin{equation}\label{phigk}
%	\left(\psi\left(\lambda_{L}\right) g_{\hat{k}-j}^{(L)}\right)^{2}>G(\hat{k}-j, L-1), ~~ j\in[1,m],
%\end{equation}
where $G(k, l)=\sum_{i=1}^{l}\left(\psi\left(\lambda_{i}\right) g_{k}^{(i)}\right)^{2}$. This together with Theorem \ref{theoremalpha} yields that 
\begin{equation*}
\begin{aligned}
\left|1-\tilde{\alpha}_{\hat{k}-1}\lambda_{L}\right| &=\tilde{\alpha}_{\hat{k}-1} \lambda_{L}-1 \leq \frac{\left(\psi(A) g_{\hat{k}-1}\right)^{T}A^2\left(\psi(A) g_{\hat{k}-1}\right)}{\left(\psi(A) g_{\hat{k}-1}\right)^{T} A^3\left(\psi(A) g_{\hat{k}-1}\right)} \lambda_{L}-1 \\ &=\frac{\sum_{i=1}^{n}\left(\lambda_{L}-\lambda_{i}\right)\lambda_{i}^2\left(\psi\left(\lambda_{i}\right) g_{\hat{k}-1}^{(i)}\right)^{2}}{\sum_{i=1}^{n} \lambda_{i}^3\left(\psi\left(\lambda_{i}\right) g_{\hat{k}-1}^{(i)}\right)^{2}} \\ & \leq \frac{\left(\lambda_{L}-\lambda_{1}\right)\lambda_{1}^2 G(\hat{k}-1), L-1)}{\lambda_{1}^3 G(\hat{k}-1), L-1)+\lambda_{L}^3\left(\psi\left(\lambda_{L}\right) g_{\hat{k}-1}^{(L)}\right)^{2}} \\ & \ \leq \frac{\lambda_{L}-\lambda_{1}}{\lambda_{L}+\lambda_{1}} < \theta. 
\end{aligned}
\end{equation*}
Using similar arguments, we can show $\left|1-\alpha_{\hat{k}}^{BB1}\lambda_{L}\right|\leq\theta$. Thus,  $\left|g_{\hat{k}}^{(L)}\right| \leq \theta\left|g_{\hat{k}-1}^{(L)}\right| \leq C_{L} \theta^{\hat{k}}$ which contradicts our assumption. Hence \eqref{Citheta} holds for all $i$. We complete the proof.
\qed
	\end{proof}

\section{Numerical experiments}\label{sec:5}
In this section, we provide numerical experiments of Algorithm \ref{Al1} on solving quadratic problems. All our codes are written in MATLAB  R2016b and  carried out on a PC with an AMD Ryzen $5$ PRO $2500$U, $2.00$ GHz processor and $8$ GB of RAM running Windows $10$ system.

For Algorithm \ref{Al1}, we use $\tilde{\alpha}_{k}$ with $\psi(A) =I$. Then, if we compute  the vector $z = Ag_{k}$ at each iteration, and keep
 into memory  $w = Ag_{k-1}$,  $\tilde{\alpha}_{k-1}$ can be expressed as
 $$c_{0}=g_{k-1}^{T} g_{k-1}, \quad c_{1}=g_{k-1}^{T} w, \quad c_{2}=w^{T} w,$$
 $$c_{3}=\frac{g_{k}^{T} z-c_{1}+2 \alpha_{k-1} c_{2}}{\alpha_{k-1}^{2}}, \quad c_{4}=\frac{z^{T} z-c_{2}+2 \alpha_{k} c_{3}}{\alpha_{k}^{2}}.$$
  Hence only  one  matrix-vector product is required in per iteration for  Algorithm \ref{Al1}.

Firstly,  we compare  Algorithm \ref{Al1}  with 
the BB1 \cite{barzilai1988two}, DY \cite{yuan2008step},  ABBmin2 \cite{frassoldati2008new}, SDC \cite{de2014efficient} and MGC \cite{zou2021fast} methods and the Alg.1 method in \cite{sun2020new} (we note this  method as SL) on solving the following quadratic problem \cite{yuan2006new}:
\begin{equation}\label{function}
	f(x)=\left(x-x^{*}\right)^{T}V\left(x-x^{*}\right),
\end{equation}
where $V = \mathrm{diag}\left\{v_{1}, \ldots, v_{n}\right\}$ is a diagonal matrix and $x^{*}$ is randomly generated with components in $[-10,10]$.  Five different
distributions of the diagonal elements $v_{j}, j=1,2, \ldots, n$, are generated,  see Table 1.

The parameter $m$ in  SL uses the value $m=5$ for the second problem set and $m = 6$ for other sets to get good performance. For ABBmin2, $\tau$ is set to $0.9$ as suggested in \cite{frassoldati2008new}. For  SDC and MGC, the parameter pairs $(h, s)$  and $(d_{1}, d_{2})$ are  all set to $(8, 6)$, which are more efficient than other choices for this test. The parameters in Algorithm \ref{Al1}  are  set to $\tau = 0.3$ and $r = 5$.

For all comparison methods, the stopping condition is $\left\|g_{k}\right\|_{2} \leq \epsilon\left\|g_{0}\right\|_{2}$,
where  $\epsilon>0$ is a given tolerance.  The problem dimension is set to $n = 1000$.  For each problem set, three different    tolerance parameters $\epsilon =  10^{-6}, 10^{-9}, 10^{-12}$ and condition numbers $\kappa = 10^{4}$, $10^{5}, 10^{6}$  are tested. For each value of $\kappa$ or $\epsilon$, 10 starting points are randomly generated   in $[-10, 10]$ and the average results  are presented in Table 2.

We see that  Algorithm \ref{Al1} clearly  outperforms the BB1, DY, SDC and MGC   methods. As compared with the SL method, Algorithm \ref{Al1} is usually faster than it when a high accuracy is required. Moreover, Algorithm \ref{Al1} often performs better than ABBmin2 for the second to fourth problem sets and is very competitive with ABBmin2 for the first and last problem sets.

\newcommand{\tabincell}[2]{\begin{tabular}{@{}#1@{}}#2\end{tabular}}
\tabcolsep12.2pt
\arrayrulewidth1pt
	\begin{table}
		\begin{center}
		\caption{Distributions of $v_{j}$}
		\label{tab:1}   
		% Give a unique label
		
		\begin{tabular}{|c|c|}
			\hline
			Problem & \  Spectrum   \\
			\hline
			1 &  \tabincell{l}{$v_{1} = 1,~v_{n} = \kappa$\\
				$\left\{v_{2}, \ldots, v_{n-1}\right\} \subset(1, \kappa)$ 
			} \\
			\hline
			2 & \tabincell{l}{$v_{1} = 1,~v_{n} = \kappa$\\$\left\{v_{2}, \ldots, v_{n / 5}\right\} \subset(1,100)$\\
				$\left\{v_{n / 5+1}, \ldots, v_{n-1}\right\} \subset\left(\frac{\kappa}{2}, \kappa\right)$
			} \\
			\hline
		
			3 & \tabincell{l}{$v_{1} = 1,~v_{n} = \kappa$\\$\left\{v_{2}, \ldots, v_{n / 5}\right\} \subset(1,100)$\\
				$\left\{v_{n / 5+1}, \ldots, v_{4 n / 5}\right\} \subset\left(100, \frac{\kappa}{2}\right)$\\
				$\left\{v_{4 n / 5+1}, \ldots, v_{n-1}\right\} \subset\left(\frac{\kappa}{2}, \kappa\right)$}  \\
			\hline
			4 & $v_{j}=\kappa^{\frac{n-j}{n-1}}$,~~
				 $j=1, \ldots, n$\\
			\hline
			
			5 & $v_{j}=\frac{\kappa}{2}\left(\cos \frac{n-j}{n-1} \pi+1\right)$,~
					~ $j=1, \ldots, n$
			 \\
		
			\hline
		\end{tabular}
	\end{center}
	\end{table}

\begin{table}
	\caption{The average number of iterations  required by Algorithm \ref{Al1}, the BB1, DY, SL, SDC, MGC and ABBmin2  methods on  quadratic problem \eqref{function} with spectral distributions  in Table \ref{tab:1}
}
	\label{tab:2}       % Give a unique label
	\footnotesize
\begin{center}
	
\tabcolsep3.6pt
\arrayrulewidth1pt
		\begin{tabular}{*{9}{|c}|}
		\hline
		$\kappa$    &      $\epsilon$  &   BB1   & DY  & SL     & MGC   & SDC     & ABBmin2& Algorithm 1 \\ \hline
		\multicolumn{9}{|c|}{Problem set 1}\\
		\hline
		
		\multirow{3}{10pt}[3pt]{$10^{4}$}  & $10^{-6}$& 323.9                & 308.1                & 249.3                & 269.4                & 262.2                & \textbf{219.8}                & 234.9                \\

		& $10^{-9}$ &	894.8                & 848.7                & 660.3                & 733.5               & 754.1                & \textbf{366.1}                & 445.8                \\

		& $10^{-12}$ &		1359.1               & 1328.6               & 1029.3               & 1142.7               & 1176.0               & \textbf{528.6}                & 649.2                \\

		\hline	
		
		\multirow{3}{10pt}[3pt]{$10^{5}$}  & $10^{-6}$&		231.5                & 238.2                & 199.7                & 204.1                & 200.3                & 203.7                & \textbf{199.3}                \\

		& $10^{-9}$ &	2265.1               & 2415.5               & 1586.9               & 1733.5               & 1806.7               & \textbf{525.6}                & 600.3                \\

		& $10^{-12}$ &		4073.0               & 4335.5               & 2786.5               & 3130.5               & 3268.0               & \textbf{687.5}                & 814.6                \\

		\hline
		
		\multirow{3}{10pt}[3pt]{$10^{6}$}  & $10^{-6}$ &
					221.0                & 222.1                & 185.3                & 196.6                & 184.3                & \textbf{173.8}                & 182.9                \\

		& $10^{-9}$ &	5221.4               & 6446.6               & 3166.5               & 3891.2               & 3800.8               & \textbf{798.0}                & 869.8                \\

		& $10^{-12}$ &	
	11834.4              & 15815.0              & 7378.4               & 8560.6               & 8276.9               & \textbf{1025.8}               & 1151.6               \\
		
		\hline 
		\multicolumn{9}{|c|}{Problem set 2}\\
		\hline
		\multirow{3}{10pt}[3pt]{$10^{4}$} & $10^{-6}$ &286.2                & 362.4                & 199.9                & 264.3                & 260.9                & 264.7                & \textbf{170.0}                \\

		& $10^{-9}$ &	653.1                & 849.7                & 430.5                & 555.9                & 556.0                & 529.0                & \textbf{350.2}                \\

		& $10^{-12}$ &		978.1                & 1301.3               & 651.2                & 840.2                & 810.0                & 786.7                & \textbf{523.9}                \\

		\hline	
		
		\multirow{3}{10pt}[3pt]{$10^{5}$}  & $10^{-6}$&		432.5                & 494.0                & 182.7                & 244.6                & 261.7                & 453.5                & \textbf{141.1}                \\

		& $10^{-9}$ &	1444.2               & 1985.0               & 616.1                & 944.2               & 900.9                & 1141.2               & \textbf{402.9}                \\

		& $10^{-12}$ &		2327.8               & 3534.1               & 1003.9               & 1427.1               & 1401.9               & 1809.0               & \textbf{635.7}                \\

		\hline
		
		\multirow{3}{10pt}[3pt]{$10^{6}$}  & $10^{-6}$ &	514.2                & 459.1                & 104.8                & 130.5                & 126.5                & 465.4                & \textbf{76.0}                 \\

		& $10^{-9}$ &3256.7               & 3427.3               & 733.6                & 925.7               & 946.0                & 2556.5               & \textbf{410.9}                \\

		& $10^{-12}$ &	
			5820.9               & 6938.1               & 1327.8               & 1790.8               & 1801.3               & 4322.4               & \textbf{742.3}                \\
		
		\hline 
		\multicolumn{9}{|c|}{Problem set 3}\\
		\hline
		\multirow{3}{10pt}[3pt]{$10^{4}$} & $10^{-6}$ &461.2                & 452.9                & \textbf{375.9}                & 404.3                & 396.5                & 401.3                & 391.6              \\

		& $10^{-9}$ &		1021.5               & 1029.5               & 821.8                &926.3               & 904.4                & \textbf{786.1}                & 790.8                \\

		& $10^{-12}$ &	1550.0               & 1550.7               & 1237.5               & 1371.1               & 1370.5               & \textbf{1142.1}               & 1165.4               \\

		\hline	
		
		\multirow{3}{10pt}[3pt]{$10^{5}$}  & $10^{-6}$&		809.7                & 804.1                & \textbf{657.4}                & 750.6                & 710.9                & 938.3                & 677.2                \\

		& $10^{-9}$ &		3089.1               & 3058.8               & 2285.1               & 2700.0               & 2555.1               & 2145.4               & \textbf{2052.6}               \\

		& $10^{-12}$ &		5128.4               & 5104.9               & 3661.4               & 4092.1               & 4014.5               & 3300.3               & \textbf{3209.4}               \\

		\hline
		
		\multirow{3}{10pt}[3pt]{$10^{6}$}  & $10^{-6}$ &		1155.6               & 1119.2               & 856.7                & 1025.7               & 982.8                & 1266.0               & \textbf{694.5}                \\

		& $10^{-9}$ &	8562.1               & 9690.6               & 6101.3               & 6876.4              & 6722.9               & 5968.5               & \textbf{4084.4}               \\

		& $10^{-12}$ &	16239.6              & 19489.3              & 10609.8              & 12340.0              & 12342.7              & 9598.3               & \textbf{6956.4}               \\

		\hline 
		\multicolumn{9}{|c|}{Problem set 4}\\
		\hline
		\multirow{3}{10pt}[3pt]{$10^{4}$} & $10^{-6}$ &628.7                & 602.1                & \textbf{497.5}                & 541.3                & 531.1                & 514.3                & 509.9                \\

		& $10^{-9}$ &		1164.8               & 1159.7               & 946.6                & 1055.1               & 1021.8               & 913.9                & \textbf{910.3}                \\

		& $10^{-12}$ &		1688.5               & 1670.4               & 1355.5               & 1507.8               & 1479.2               & 1296.5               & \textbf{1291.3}               \\

		\hline	
		
		\multirow{3}{10pt}[3pt]{$10^{5}$}  & $10^{-6}$&			1480.7               & 1442.2               & 1118.1               & 1304.5               & 1297.6               & 1330.6               & \textbf{1177.2}               \\

		& $10^{-9}$ &		3440.9               & 3702.2               & 2705.9               & 3056.5              & 3038.8               & 2664.7               & \textbf{2469.3}               \\

		& $10^{-12}$ &			5437.8               & 5795.1               & 4240.0               & 4665.7               & 4582.6               & 3891.6               & \textbf{3718.5}               \\

		\hline
		\multirow{3}{10pt}[3pt]{$10^{6}$}  & $10^{-6}$ &		2533.3               & 2667.0               & 1941.2               & 2422.4               & 2202.7               & 3119.5               & \textbf{1996.0}               \\

		& $10^{-9}$ &	10849.6              & 12898.5              & 7542.4               & 9075.5             & 8508.5               & 7713.0               & \textbf{6570.9}               \\

		& $10^{-12}$ &	18028.8              & 23430.2              & 12681.4              & 14396.4              & 13549.1              & 11838.3              & \textbf{10539.5}              \\
		
		\hline 
		\multicolumn{9}{|c|}{Problem set 5}\\
		\hline
		\multirow{3}{10pt}[3pt]{$10^{4}$} & $10^{-6}$ &527.0                  & 535.3               & \textbf{454.3}               & 487.0                  & 487.7               & 599.6               & 484.5               \\

		& $10^{-9}$ &			4455.5              & 5480.7              & 4115.5              & 4098.2              & 3805.4              & \textbf{2879.0}              & 2993.8              \\

		& $10^{-12}$ &			9098.2              & 10185.5             & 7358.3              & 7815.6              & 7473.1              & \textbf{4312.4}              & 4806.7              \\

		\hline	
		
		\multirow{3}{10pt}[3pt]{$10^{5}$}  & $10^{-6}$&			552.7               & 573.5               & \textbf{467.7 }              & 525.8               & 497.4               & 635.6                & 489.5               \\

		& $10^{-9}$ &			4633.8              & 5142.8              & 3812.0                 & 4047.6              & 3888.3              & \textbf{2837.5}              & 2879.0              \\

		& $10^{-12}$ &			9007.1              & 10727.7             & 7141.0                 & 7150.1              & 6967.7              & \textbf{4209.8}              & 4767.2               \\

		\hline
		\multirow{3}{10pt}[3pt]{$10^{6}$}  & $10^{-6}$ &			527.5               & 558.6               & \textbf{476.9}              & 503.0                  & 477.0                  & 669.9               & 489.9               \\

		& $10^{-9}$ &	4859.5              & 5423.8              & 4121.7              & 4127.7              & 4363.2              & \textbf{2846.0}                 & 3060.0              \\

		& $10^{-12}$ &			9410.7              & 10661.0             & 7703.8              & 7299.8              & 7528.9              & \textbf{4379.8}              & 5083.6              \\
		
		\hline 
	\end{tabular}
\end{center}
\end{table}

\tabcolsep20pt
\arrayrulewidth1pt
\begin{table}
\begin{center}
	\caption{Test problems from Suitesparse Matrix Collection}
	\label{tab:4}       % Give a unique label
	\begin{tabular}{|c|c|c|c|}
		\hline
	 Matrices   &  Size  &  Nonzeros  &Condition number  \\
		 \hline
		bcsstk14&   1806  & 63454  & $1.192324\times 10^{10}$  \\
		bcsstk15& 3948  & 117816  & $6.538185\times 10^9$ \\
	   bcsstk17& 10974 & 428650   & $1.296064\times 10^{10}$\\
	    bcsstk18& 11948 & 149090  & $3.459995\times 10^{11}$\\
	   msc01440& 1440  & 44998   & $3.305875 \times 10^6$ \\
       msc04515& 4515  & 97707   & $	2.272772 \times 10^6$ \\
	    ex15& 6867  & 98671   & $8.612330 \times 10^{12}$ \\
	    cbuckle& 13681 & 676515   & $3.299134\times 10^7$\\
        
	   gyro\_k& 17361 & 1021159 & $1.095832\times 10^9$\\
	    s3dkq4m2& 90449 & 4427725 & $1.896133\times 10^{11}$\\
	    s3dkt3m2& 90449 & 3686223 & $3.625322\times 10^{11}$ \\
	   s3rmq4m1& 5489  & 262943   &$1.765559\times 10^{10}$\\
        s3rmt3m1& 5489  & 217669   & $2.481977\times 10^{10}$\\
       s3rmt3m3& 5357  & 207123   & 	$2.400640 \times 10^{10}$\\
                      
		\hline
	\end{tabular}
\end{center}
\end{table}

\tabcolsep5pt
\arrayrulewidth1pt

\begin{table}
	\begin{center}
	\caption{The average number of iterations  required by Algorithm \ref{Al1}, the BB1, DY, SL, SDC, MGC and ABBmin2 methods  on problem \eqref{1a}  with $A$ given by Table \ref{tab:4}}
	\label{tab:5}       % Give a unique label
	\begin{tabular}{|c|c|c|c|c|c|c|c|}
		\hline
	    Matrices     &  BB1  &DY  &SL      &MGC    &SDC  &ABBmin2 &Algorithm 1 \\
		\hline
		
	    bcsstk14&2425.8               & 1948.7               & 1791.0                 & 1705.9               & 1820.6               & \textbf{1667.4}               & 1677.3               \\

		bcsstk15&	    2510.4               & 2219.6               & 1966.0                 & 1778.3               & 2000.3               & 1788.1               & \textbf{1774.2}               \\

		bcsstk17& 		4411.5               & 4335.5               & 3721.6               & 3428.0                 & 3589.0                 & 3035.4               & \textbf{3005.4}               \\

		bcsstk18& 			4118.8               & 4194.9               & 3498.3               & 3655.1               & 4055.5               & 3334.6               & \textbf{3206.7}               \\

	msc01440&			1423.8               & 1623.9               & 1329.4               & 1233.0                 & 1241.5               & 1151.9               & \textbf{1139.0}                 \\

		msc04515& 	1008.8               & 986.2                & 903.7                & \textbf{810.7}                & 851.6                & 821.4                & 829.7                \\

		ex15   &			913.0                  & 1002.5               & 802.6                & 806.7                & 778.9                & \textbf{716.9}                & 775.4                \\

		cbuckle&		2335.5               & 2092.8               & 2015.9               & 1884.0                 & 1984.6               & \textbf{1665.4}               & 1730.0                 \\

		gyro\_k& 			1672.8               & 1564.6               & 1461.1               & 1390.7               & 1366.9               & \textbf{1321}                 & 1329.2               \\

		s3dkq4m2&			1611.8               & 1800.6               & 1600.3               & 1362.8               & 1493.4               & 1205.1               & \textbf{1178}                 \\

		s3dkt3m2& 		2048.6               & 2081.8               & 1802.5               & 1455.3               & 1644.5               & \textbf{1158.6}               & 1236.5               \\

		s3rmq4m1& 			1153.7               & 1191.6               & 967.1                & 880.6                & 911.9                & 835.7                & \textbf{811.5}                \\

		s3rmt3m1&			1038.7               & 1099.9               & 929.4                & 871.6                & 906.4                & \textbf{879.8}                & 904.1                \\

		s3rmt3m3& 			1606.7               & 1612.2               & 1457.9               & 1276.6               & 1286.6               & 1229.6               & \textbf{1206.6}               \\
		              
	\hline
	\end{tabular}
\end{center}
\end{table}

\tabcolsep6.4pt
\arrayrulewidth1pt

\begin{table}
\begin{center}	 
	\caption{The CPU time  required by Algorithm \ref{Al1}, the  BB1, DY, SL, SDC, MGC and ABBmin2 methods on problem \eqref{1a} with $A$ given by Table \ref{tab:4}}
	\label{tab:5}       % Give a unique label
	\begin{tabular}{|c|c|c|c|c|c|c|c|}
		\hline
		Matrices     &  BB1  &DY  &SL      &MGC    &SDC  &ABBmin2 & Algorithm 1  \\
		\hline
		bcsstk14&0.40                 & 0.17                 & 0.16                 & 0.31                 & 0.16                 & 0.30                 & \textbf{0.15}                 \\

		bcsstk15& 	0.57                 & 0.28                 & 0.25                 & 0.41                 & 0.24                 & 0.43                 & \textbf{0.22}                 \\

		bcsstk17& 		6.44                 & 3.51                 & 2.99                 & 5.04                 & 2.81                 & 4.49                 & \textbf{2.41}                 \\

		bcsstk18& 		3.27                 & 2.04                 & 1.69                 & 3.02                 & 1.79                 & 2.80                 & \textbf{1.54}                 \\

		msc01440&	0.13                 & 0.09                 & 0.07                 & 0.12                 & \textbf{0.06}                 & 0.11                 & \textbf{0.06}                 \\

		msc04515&		0.21                 & 0.12                 & 0.11                 & 0.19                 & \textbf{0.10}                 & 0.18                 & \textbf{0.10}                \\

		ex15   &		0.20                 & 0.13                 & \textbf{0.10}                 & 0.19                 & \textbf{0.10}                 & 0.17                 & \textbf{0.10}                 \\

		cbuckle&		6.12                 & 2.98                 & 2.87                 & 5.00                 & 2.67                 & 4.32                 & \textbf{2.29}                 \\

		gyro\_k&		5.26                 & 2.65                 & 2.49                 & 4.44                 &\textbf{2.27}                 & 4.23                 & \textbf{2.27}                 \\

		s3dkq4m2& 		26.91                & 16.35                & 14.33                & 23.59                & 12.91                & 20.97                & \textbf{10.60}                \\

		s3dkt3m2& 		25.26                & 13.87                & 12.39                & 19.27                & 10.85                & 15.20                & \textbf{8.62}                 \\

		s3rmq4m1& 		0.60                 & 0.32                 & 0.26                 & 0.51                 & 0.25                 & 0.51                 & \textbf{0.22}                 \\

		s3rmt3m1& 		0.48                 & 0.27                 & 0.28                 & 0.43                 & \textbf{0.22}                & 0.42                 & \textbf{0.22}                 \\

		s3rmt3m3& 			0.79                 & 0.44                 & 0.39                 & 0.67                 & 0.33                 & 0.62                 & \textbf{0.32}                 \\

		\hline
	\end{tabular}
\end{center}
\end{table}

Next, we compare the  above methods on problem \eqref{1a}, where $b = Ax^{*}$ and  $x^{*}$ is a random vector as before. We test 14 sparse matrices from the SuiteSparse Matrix Collection \cite{davis2011university}
listed in Table 3. The iteration stops when  $\left\|g_{k}\right\|_{2} \leq 10^{-6}\left\|g_{0}\right\|_{2}$. In our test, we choose the parameters so that each method achieves the best performance. Specifically, we set $m = 4$ for SL,  $(h, s) = (3, 4)$ for SDC, $(d_{1},d_{2}) = (4,4)$ for MGC, $\tau = 0.3$ for ABBmin2 and  $\tau  = 0.1~\text{ and }~r = 5$ for Algorithm \ref{Al1} .

For each matrix, ten initial points between $-10$ and $10$ are randomly generated. The average number of iterations of compared methods are listed in Table 4. We see that  for most matrices  Algorithm \ref{Al1}  has better performance than  BB1, DY, SL, MGC, SDC and  is competitive with ABBmin2 in the sense of number of iterations. Table 5 lists the average CPU time in seconds for those methods. We observe that  Algorithm \ref{Al1}  takes less CPU time than  BB1, DY, MGC, ABBmin2, and is as fast as SL and SDC. One important reason for this phenomenon is that our method reuses the  stepsize $\tilde{\alpha}_{k-1}$ for $r$ iterations which can reduce computational cost.

\section{Conclusions and discussions}\label{sec:6}
For a uniform analysis, we  considered a family of gradient methods whose stepsize is provided by $\alpha_{k}$ in \eqref{alphaAOfamily}, which includes the Dai--Yang method \eqref{dystepsize} as a special case. It is proved the family zigzags between two directions in a subspace spanned by the two eigenvectors corresponding to  the smallest and largest eigenvalues of the Hessian. In order to achieve the  two dimensional quadratic termination of the family, we derived a short stepsize $\tilde{\alpha}_{k}$ \eqref{alpha1} that converges to $1/\lambda_{n}$.  By  using the long BB stepsize and  $\tilde{\alpha}_{k-1}$  in a new adaptive cyclic manner, we designed Algorithm \ref{Al1}  for unconstrained quadratic optimization. We proved that Algorithm \ref{Al1} converges $R$-linearly at a rate of $1-1/\kappa$. Our numerical results  on minimizing quadratic functions indicate  the efficiency of Algorithm \ref{Al1} over other recent successful gradient methods.

%We  presented the Property B and proved that any gradient method satisfying this property has $R$-linear convergence which implies Algorithm \ref{Al1} and many others are $R$-linearly convergent. Under  proper conditions, we proved the  $R$-factor  can be bounded above by $1-1/\kappa$. 

By using the same arguments as those in the proof of  Theorem \ref{Finite termination}, we find  stepsizes $\tilde{\beta}_{k}$ and $\tilde{\gamma}_{k} $ are such that \eqref{Falpha1} and \eqref{Falpha2} achieve the
 two dimensional quadratic termination, respectively. For the $n$ dimensional quadratic problem, by  Theorem 1 in \cite{huang2019asymptotic}, we obtain 
$$\lim _{k \rightarrow \infty} \tilde{\beta}_{k}=\frac{1}{\lambda_{n}},~\lim _{k \rightarrow \infty} \tilde{\gamma}_{k}=\frac{1}{\lambda_{n}}$$
and 
$$\lim _{k \rightarrow \infty} \hat{\beta}_{k}=\frac{1}{\lambda_{1}},~\lim _{k \rightarrow \infty} \hat{\gamma}_{k}=\frac{1}{\lambda_{1}}.$$ 
In addition, we have $\lim _{k \rightarrow \infty} \hat{\alpha}_{k}=1/\lambda_{1}$. However, $\hat{\alpha}_{k}$, $\hat{\beta}_{k}$ and $\hat{\gamma}_{k}$ would not be  good approximations of $1/\lambda_{1}$, for more details see \cite{huang2020gradient}.  It is worth noting that   the  stepsize  proposed in \cite{frassoldati2008new} is a special case of $\tilde{\beta}_{k}$ with $\psi(A)=I$. Moreover, it is not difficult to prove that  $\tilde{\beta}_{k}$ and $\tilde{\gamma}_{k}$
are short stepsizes in the sense $\tilde{\beta}_{k}<c_{2}/c_{3}$ and $\tilde{\gamma}_{k}<c_{1}/c_{2}$.  Hence, we can replace $\tilde{\alpha}_{k-1}$ in Algorithm \ref{Al1} by $\tilde{\beta}_{k-1}$ and $\tilde{\gamma}_{k-1}$, which leads to  two variants of Algorithm \ref{Al1}. Preliminary experimental results show that the two variants are competitive with Algorithm \ref{Al1}.  Furthermore, we can obtain  the same convergence results as Theorem \ref{converth2}.

The results of this paper show that the two dimensional quadratic termination property is useful for designing efficient gradient methods. It would be interesting to develop new algorithms for solving general unconstrained optimization problems based on such a property. We leave this as our future work.

\section{Declarations}
All data generated or analyzed during this study are included in this published article and
are also available from the corresponding author on reasonable request.

\begin{acknowledgements}
This work was supported by the National Natural Science Foundation of China  (Grant No. 11701137) and Natural Science Foundation of Hebei Province (Grant No. A2021202010).
\end{acknowledgements}


%\newpage
\bibliographystyle{spmpsci}      % mathematics and physical sciences

%\begin{thebibliography}{}
% etc
%\end{thebibliography}
%\bibliography{BB_hyk}

\end{document}